\documentclass[12pt,a4paper]{amsart}
\usepackage{amssymb,amsmath,amsthm}
\usepackage{graphicx}
\usepackage{hyperref}
 \usepackage[color=green!40]{todonotes}

\newcommand\cA{{\mathcal A}}

\newcommand\cF{{\mathcal F}}
\newcommand\cS{{\mathcal S}}
\newcommand\x{{\mathbf x}}
\newcommand\y{{\mathbf y}}
\newcommand\z{{\mathbf z}}

\makeatletter
\newtheorem*{rep@theorem}{\rep@title}
\newcommand{\newreptheorem}[2]{%
\newenvironment{rep#1}[1]{%
 \def\rep@title{#2 \ref{##1}}%
 \begin{rep@theorem}}%
 {\end{rep@theorem}}}
\makeatother

\theoremstyle{plain}
\newtheorem{theorem}{Theorem}[section]
\newreptheorem{theorem}{Theorem}
\newtheorem{lemma}[theorem]{Lemma}
\newtheorem{corollary}[theorem]{Corollary}

\newtheorem{proposition}[theorem]{Proposition}

\newtheorem{claim}[theorem]{Claim}
\theoremstyle{definition}

\newtheorem{defn}[theorem]{Definition}

\newcommand\cref[1]{Corollary~\ref{cor:#1}}

\newcommand{\ex}{\mathop{}\!\mathrm{ex}}

\title{Extremal graph theoretic questions for q-ary vectors}

\author{Bal\'azs Patk\'os}
\address{Alfr\'ed R\'enyi Institute of Mathematics}
\email{patkos@renyi.hu}
\thanks{Patk\'os's research is partially supported by NKFIH grants SNN 129364 and FK 132060.}

\author{Zsolt Tuza}
\address{Alfr\'ed R\'enyi Institute of Mathematics and University of Pannonia}
\email{tuza.zsolt@mik.uni-pannon.hu}
\thanks{Tuza's research is partially supported by NKFIH grant SNN 129364.}

\author{M\'at\'e Vizer}
\address{Alfr\'ed R\'enyi Institute of Mathematics}
\email{vizermate@gmail.com}
\thanks{Vizer's research is partially supported by NKFIH grants SNN 129364, FK 132060, KH130371, by the  J\'anos Bolyai Research Fellowship  and by the New National Excellence Program under the grant number \'UNKP-21-5-BME-361.}

\date{}
\begin{document}

\maketitle

\begin{abstract}
   A $q$-graph $H$ on $n$ vertices is a set of vectors of length $n$ with all entries from $\{0,1,\dots,q\}$ and every vector (that we call a $q$-edge) having exactly two non-zero entries. The support of a $q$-edge $\x$ is the pair $S_\x$ of indices of non-zero entries. We say that $H$ is an $s$-copy of an ordinary graph $F$ if $|H|=|E(F)|$, $F$ is isomorphic to the graph with edge set $\{S_\x:\x\in H\}$, and whenever $v\in e,e'\in E(F)$, 
   the entries with index corresponding to $v$ in the $q$-edges corresponding to $e$ and $e'$ sum up to at least $s$. E.g., the $q$-edges $(1,3,0,0,0), (0,1,0,0,3)$, and $(3,0,0,0,1)$ form a 4-triangle.
   
   The Tur\'an number $\ex(n,F,q,s)$ is the maximum number of $q$-edges that a $q$-graph $H$ on $n$ vertices can have if it does not contain any $s$-copies of $F$. In the present paper, we determine the asymptotics of $\ex(n,F,q,q+1)$ for many graphs $F$.
\end{abstract}

\section{Introduction}

The main problem of extremal graph theory is to determine the Tur\'an number $\ex(n,F)$, the maximum number of edges that an $n$-vertex graph $H$ can have if $H$ does not contain $F$ as a subgraph. The analogous parameter for forbidding a family $\cF$ of subgraphs is denoted by $\ex(n,\cF)$. The asymptotics of $\ex(n,F)$ is given by the Erd\H os--Stone--Simonovits theorem \cite{ESt,ES} whenever $F$ has chromatic number at least 3. For many bipartite graphs $F$, even the order of magnitude is unknown; for a survey see \cite{FS}.

\begin{theorem}[Erd\H os--Stone--Simonovits, \cite{ES,ESt}]\label{ESS}
For a family\/ $\cF$ of graphs, let\/ $\chi(\cF)$ denote the minimum value of\/ $\chi(F)$ over all members\/ $F$ of\/ $\cF$. Then\/ $\ex(n,\cF)=(1-\frac{1}{\chi(\cF)-1}-o(1))\binom{n}{2}$. 
\end{theorem}

In this paper, we present a generalization of this problem via the incidence matrix of graphs. The \textit{incidence matrix} of a graph $G$ is an $|E(G)| \times |V(G)|$ matrix $M$, whose rows are indexed with the edges, the columns are indexed with the vertices of $G$, and for $v \in V(G)$ and $e \in E(G)$ the corresponding entry is $M[e,v]=1$ if $v \in e$ and $M[e,v]=0$ if $v \notin e$. A 0-1 matrix is an incidence matrix of a simple graph if and only if all its columns are pairwise distinct, and every column contains exactly two 1-entries. 

We want the setting to be more general by allowing entries from $\{0,1,\dots,q\}$, but still every column should have exactly two non-zero entries. To do so first we define the generalization of an edge on $n$ vertices: it is a vector $\x \in \{0,1,\ldots,q\}^n$ that we call a \emph{$q$-edge}, and its $i$th coordinate will be denoted by $x_i$. The \textit{support} of $\x\in \{0,1,\ldots,q\}^n$ is $S_\x:=\{i: 1 \le i \le n, \ x_i \neq 0\}$.
We introduce the notation $Q(n,r):= \{ \x \in \{0,1,\ldots,q\}^n$ : $|S_\x|=r \}$. Then a \textit{$q$-graph} on $n$ vertices is any subset $H$ of $Q(n,2)$.
For a $q$-graph $H$, we write $S_H$ to denote $\{S_\x:\x\in H\}$. 
We denote $Q(n,2)$ by $2(n,2), 3(n,2)$ in case $q=2,3$ respectively.

We still need to define when two $q$-edges intersect at a certain vertex, or more generally when a $q$-graph $H$ contains an ordinary graph $F$. There are probably several ways how one could do that, we use the following approach of ours: in the recent manuscript \cite{PTV}, generalizing set intersection, for an integer $s \ge 1$ we defined the \textit{$s$-sum intersection} of two $q$-edges $\x,\y \in \{0,1,\ldots,q\}^n$ as $\x \cap_s \y = \{i: 1 \le i \le n \textrm{ for which } x_i + y_i \ge s\}$. (The case $q=1$, $s=2$ is ordinary set intersection via characteristic vectors.) As an analogous definition, we would say that $q$-edges of support size 2 form an $s$-copy of a graph $F$ if they $s$-sum intersect ``according to $F$''. Formally, we introduce the following definition.

\begin{defn}

If\/ $F=(V(F),E(F))$ is an ordinary graph without isolated vertices and\/ $H \subseteq Q(n,2)$ is a\/ $q$-graph, then we say that\/ $H$ is an\/ \textit{$s$-copy of $F$} if the graph\/ $(V,E)$ with\/ $V=\cup_{\x\in H}S_\x$, $E=S_H$ is isomorphic to\/ $F$; and if\/ $\iota:F\rightarrow (V,E)$ is an isomorphism\footnote{We call $\iota:F\rightarrow (V,E)$ an isomorphism if $\iota: V(F) \rightarrow V$ is a bijection that induces $\iota: E(F) \rightarrow E$ by $\iota(e):=\{ \iota(v) : v \in e\}$ for $e \in E(F)$.}, then for every\/ $e,e'\in E(F)$ and\/ $v\in e\cap e'$ we have\/ $\iota(e)_{\iota(v)}+\iota(e')_{\iota(v)}\ge s$.

For graphs\/ $F$ with isolated vertices, we say that\/ $H\subseteq Q(n,2)$ is an\/ \textit{$s$-copy of $F$} if\/ $n\ge |V(F)|$ and\/ $H$ is an\/ $s$-copy of\/ $F[U]$, where\/ $U$ is the set of non-isolated vertices of\/ $F$.


\end{defn}

For example the $q$-edges $\x,\y,\z \in 3(n,2)$ with $x_1=y_2=z_3=1$ and $x_2=y_3=z_1=3$ form a 4-copy of a triangle, but $\x',\y',\z' \in 3(n,2)$ with $x'_2=y'_2=z'_3=1$ and $x'_1=y'_3=z'_1=3$ do not since $x'_2+y'_2=2<4$.

\begin{defn} For a(n ordinary) graph\/ $F$ and\/ $n,q,s \ge 1$ we introduce
\[
\ex(n,F,q,s):= \max \{|H|: H \subset Q(n,2), H \textrm{ does not contain an\/ $s$-copy of\/ }F\}.
\]

In this paper, we will only consider cases\/ $s=q+1$ and for the sake of simplicity we write\/ $\ex(n,F,q)$ for\/ $\ex(n,F,q,q+1)$.
\end{defn} 

There are several very natural first thoughts that one can have when starts thinking about determining $\ex(n,F,q)$. 

A trivial construction is to consider an ordinary $F$-free graph $G$ on $n$ vertices with $\ex(n,F)$ edges and define $H_{G,q}:=\{\x\in Q(n,2):S_\x\in E(G)\}$. Clearly, $H_{G,q}$ does not contain any $(q+1)$-copies of $F$ and $|H_{G,q}|=q^2\cdot \ex(n,F)$. This shows that $\ex(n,F,q)\ge q^2\cdot \ex(n,F)$ for any graph $F$ and $q \ge 1$, but it turns out that this inequality is hardly ever sharp. Trivial exceptions are matchings for which the  ``sum-intersection condition'' is void.

As we will show, there is a simple construction and an easy observation that partition graphs according to the asymptotics of $\ex(n,F,q)$. 

Our first result shows that the  construction $U_{q,n}$ separates the asymptotics of $\ex(n,F,q)$ for $F$ being a forest from the value of $\ex(n,F,q)$ for all graphs containing a cycle We defer the definition of $U_{q,n}$ to the next section.

\begin{theorem}\label{gentree}
For any tree\/ $T$ and integer\/ $q\ge 2$, there exists a positive constant\/ $c=c_{T,q}$ such that\/ $$\ex(n,T,q)\le |U_{q,n}|-c\binom{n}{2}=\left(\lfloor \frac{q^2}{2}\rfloor-c\right)\binom{n}{2}.$$
\end{theorem}

Our next theorem shows that the asymptotics of $\ex(n,T,2)$ (and probably of $\ex(n,T,q)$ in general, too) depends on the radius \footnote{The \textit{radius} of a connected graph $G$ is $\min_{v \in V(G)} \{ \max_{u \in V(G)} distance (u,v)\}$.} of the tree. Let $t_{n,r}$ denote the number of edges of the $r$-partite Tur\'an graph \footnote{The \textit{Tur\'an graph} $T(n,r)$ is the $n$-vertex complete $r$-partite graph whose vertex partition is as balanced as possible.} on $n$ vertices, and let $t'_{n,r}$ denote the number of edges of the complete $r$-partite graph on $n$ vertices with one part having size $\lfloor \frac{n}{2r-1}\rfloor$ with the other part sizes differing by at most one (so their size is $(1+o(1))\frac{2n}{2r-1}$).

\begin{theorem}\label{2tree}
Suppose\/ $T$ is a tree of radius\/ $r$. 
\begin{enumerate}
    \item 
    If the diameter of\/ $T$ is $2r$, then\/ $\ex(n,T,2)=(1+o(1))(\binom{n}{2}+t_{n,r})$.
    \item
    If the diameter of\/ $T$ is\/ $2r-1$, then\/ $\ex(n,T,2)=(1+o(1))(\binom{n}{2}+t'_{n,r}-\binom{\lfloor\frac{n}{2r-1}\rfloor}{2})$.
\end{enumerate}
\end{theorem}

Our next result states that the construction $U_{q,n}$ is asymptotically optimal for bipartite graphs with all components containing at most one cycle.

\begin{theorem}
\label{bipuni}
Suppose\/ $F$ is a bipartite graph such that all its components are unicyclic\footnote{A \textit{unicyclic} graph is a connected graph containing exactly one cycle.} or trees, and\/ $F$ contains at least one cycle. Then\/ $$\ex(n,F,q)=(\lfloor \frac{q^2}{2}\rfloor+o(1))\binom{n}{2}=|U_{q,n}|+o(n^2).$$
\end{theorem}

For graphs $F$ that contain an odd cycle, we have $\ex(n,F,q)\ge q^2\cdot \ex(n,F)\ge (\frac{q^2}{2}+o(1))\binom{n}{2}$. If $q$ is odd, then the coefficient is strictly larger than that of the size of $U_{q,n}$, while if $q$ is even, then the two asymptotics are the same. We show that for the triangle, this trivial lower bound is sharp for $q=2$.

\begin{theorem}\label{triangle}
For\/ $n\ge 2$ we have\/ $$\ex(n,C_3,2)=4\cdot \ex(n,C_3)=4\cdot \lfloor \frac{n^2}{4}\rfloor.$$
\end{theorem}

Whenever the forbidden graph $F$ contains a component with more than one cycle, or equivalently with more edges than vertices, then the asymptotics of $\ex(n,F,q)$ changes dramatically. This follows from a simple observation. We need the following definition: for any $H \subseteq Q(n,2)$, we denote by $H^L$ the $q$-graph formed by the $q$-edges $\x\in H$ with $x_i\ge \frac{q+1}{2}$ for all $i\in S_\x$. 

\begin{proposition}\label{nagyel}
Suppose\/ $F$ contains a component\/ $C$ such that $|V(C)|<|E(C)|$. Then a\/ $(q+1)$-copy of\/ $F$ must contain an\/ $\x\in Q(n,2)^L$. In particular,\/  $Q(n,2)\setminus Q(n,2)^L$ does not contain a\/ $(q+1)$-copy of\/ $G$ and so\/ $\ex(n,F,q)\ge (q^2-(q-\lceil \frac{q+1}{2}\rceil +1)^2)\binom{n}{2}$.
\end{proposition}

A trivial upper bound using $H^L$ is as follows. 

\begin{proposition}\label{upperL}
For any graph\/ $F$, if\/ $H\subseteq Q(n,2)$ does not contain a\/ $(q+1)$-copy of\/ $F$, then\/ $$|H|\le (q^2-(q-\lceil \frac{q+1}{2}\rceil +1)^2)\binom{n}{2}+(q-\lceil \frac{q+1}{2}\rceil +1)^2\ex(n,F).$$
\end{proposition}

\begin{proof}
This follows from the simple observation that $S_{H^L}$ is $F$-free.
\end{proof}

Proposition \ref{upperL} and Theorem \ref{ESS} yield a general upper bound involving $\chi(F)$. In the next section, we introduce the \textit{robust chromatic number} $\chi_1(F)$ of $F$ and obtain lower bounds using this parameter. We will determine the asymptotics of $ex(n,F,q)$ whenever $\chi_1(F)=\chi(F)$. 
 We will show that in some sense, almost all graphs have this property.

Let $K(m,r,p)$ denote the probability space of all labeled $r$-partite graphs with each partite set having size $m$, where any two vertices in different parts are joined with probability $p$ independently of any other pairs.

\begin{theorem}\label{random}
If\/ $p=\omega(m^{-1/\binom{r}{2}})$, then\/ $\chi_1(K(m,r,p))=r=\chi(K(m,r,p))$ with probability tending to one as\/ $m$ tends to infinity.
\end{theorem}

 One might wonder whether out of the two expressions involving the robust and the ordinary chromatic number, one always gives the asymptotics of $\ex(n,F,q)$. 
By considering complete 3-partite graphs
 we show that this is not the case. We determine the asymptotics of $\ex(n,K_{r,s,t},2)$ and will compare the next theorem to $\chi_1(K_{r,s,t})$ and $\chi(K_{r,s,t})$ in the next section. Note that the case $r=s=t=1$ was included in Theorem \ref{triangle}.

\begin{theorem}   \label{t:tripart}
For integers\/ $1\le r\le s\le t$ with\/ $t\ge 2$, we have the following:
\begin{enumerate}
    \item 
    if\/ $r=1$ or\/ $s\le 2$, then\/ $\ex(n,K_{r,s,t},2)=(3+o(1))\binom{n}{2}$,
    \item
    if\/ $r=2$ and\/ $s\ge 3$, then\/ $\ex(n,K_{2,s,t},2)=(\frac{13}{4}+o(1))\binom{n}{2}$,
    \item
    if\/ $r\ge 3$, then\/ $\ex(n,K_{r,s,t},2)=(\frac{7}{2}+o(1))\binom{n}{2}$.
\end{enumerate}
\end{theorem}

For later reference, let us state a classical result of K\H ov\'ari, S\'os, and Tur\'an that will be often used in our proofs.

\begin{theorem}[K\H ov\'ari, S\'os, Tur\'an \cite{KST}]\label{kst}
For any fixed\/ $s\le t$, we have\/ $\ex(n,K_{s,t})=O(n^{2-1/s})$.
\end{theorem}

\noindent \textbf{Notation.} In many of the proofs, for a $q$-graph $H\subseteq Q(n,2)$, we will consider the graph of pairs $(ij)$ that are supports of some edges of $H$ with prescribed entries. Formally,
for $H \subseteq Q(n,2)$ and $(a,b)\in [q]^2$, let $\overrightarrow{H}_{a,b}$ be the directed graph of all edges $(ij)$ for which the $q$-edge $\x\in Q(n,2)$ with $S_\x=\{i,j\}$ and $x_i=a,\ x_j=b$ belongs to $H$. Observe that $\overrightarrow{H}_{a,b}$ might contain directed cycles of length 2. Also, $|H|=\sum_{1\le a< b\le q}|\overrightarrow{H}_{a,b}|+\sum_{1\le a\le q}|H_{a,a}|$ will be often used to obtain bounds on $|H|$. Sometimes we will consider pairs that are supports for two prescribed pairs of entries. Formally, we define $\overrightarrow{H}_{(a,b),(c,d)}$ to be $\overrightarrow{H}_{a,b}\cap \overrightarrow{H}_{c,d}$. The pair $(c,d)$ will always be either $(b,a)$ or $(\overline{a},\overline{b})$, where $\overline{\alpha}:=q+1-\alpha$ for any $1\le \alpha\le q$.

 Finally, let $H_{a,b}$ and $H_{(a,b),(c,d)}$ denote the simple graph obtained from $\overrightarrow{H}_{a,b}$ and $\overrightarrow{H}_{(a,b),(c,d)}$ by removing orientations and multiple edges. So $\{i,j\}\in H_{(a,b),(c,d)}$ if and only if either $\x,\y\in H$ with $x_i=a,\ y_i=c,\ x_j=b,\ y_j=d$ or $\x',\y'\in H$ with $x'_i=b,\ y'_i=d,\ x'_j=a,\ y'_j=c$.

\section{Constructions, preliminary results}

We start with the construction of the universal tree $U_{q,n}$, which is very similar to the one showing that extremal problems for  oriented graphs are only meaningful if the forbidden subgraph is a directed acyclic graph: if the vertex set is $[n]=\{1,2,\dots,n\}$ and we orient every edge $(uv)$ of the complete graph on $[n]$ towards $v$ if $u<v$, then this oriented complete graph does not contain any directed cycles.

\begin{defn}\label{utre}
For\/ $n,q \ge 1$ we define the universal\/ $q$-tree\/ $U_{q,n}=U^<\cup\bigcup_{1\le i<j\le n}U^{i,j,=}$, where\/
\[
U^<=\{\x\in Q(n,2):\sum_ix_i<q+1\}, \ and\/ 
\]
\[
U^{i,j,=}=\{\x\in Q(n,2):S_\x=\{i,j\}, \ x_i< x_j,\ x_i+x_j=q+1\}.
\]
\end{defn}

\begin{lemma}
The universal\/ $q$-tree does not contain a\/ $(q+1)$-copy of any cycle.
\end{lemma}

\begin{proof}
If $C$ is a $(q+1)$-copy of a cycle, then for any $\x\in C$, we must have $\sum_{h=1}^nx_h=q+1$ as for any vertex $v$ of the cycle, the two adjacent $q$-edges $\x,\y$ we must have $x_v+y_v\ge q+1$ so $\sum_{\x\in C}\sum_{h=1}^nx_h=\sum_{v\in S_C}\sum_{\x\in C}x_v\ge |S_C|(q+1)=|C|(q+1)$. 

But if $C$ consists only of $q$-edges $\x$ with $\sum_hx_h=q+1$, then if $v$ denotes the smallest vertex in $S_C$, then for both adjacent $q$-edges $\x,\y \in C$ we should have $x_v,y_v<\frac{q+1}{2}$ contradicting $x_v+y_v\ge q+1$.
\end{proof}

We continue with the simple lower bound of Proposition \ref{nagyel}.

\begin{proof}[Proof of Proposition \ref{nagyel}]
 Without loss of generality we can assume that $F$ is connected. For any $i$ and a $q$-copy $F_q$ of $F$ there can be at most one $\x\in F_q$ with $x_i<\frac{q+1}{2}$. Therefore, the total number of $\x\in F_q$ with at least one non-zero coordinate being less than $\frac{q+1}{2}$ is at most $|V(F)|<|E(F)|$, and thus at least one $\x\in F_q$ must have both non-zero coordinates at least $\frac{q+1}{2}$. 
\end{proof}

Proposition \ref{nagyel} can be strengthened via the following definition. 

\begin{defn}

We say that a mapping\/ $f:V(F)\rightarrow E(F)$ is a\/ \textit{1-selection} of the graph\/ $F$ if for every vertex\/ $v\in V(F)$ we have\/ $v\in f(v)$. The\/ \textit{1-removed} $F_f$ of $F$ with respect to the 1-selection\/ $f$ is the graph with\/ $V(F_f)=V(F)$ and\/ $E(F_f)=E(F)\setminus f[V(F)]$. Let\/ $\cF_F$ denote the set of all 1-removed graphs of\/ $F$.
The\/ \textit{$1$-robust chromatic number} $\chi_1(F)$ is\/ $\min\{\chi(F_f): F_f\in \cF_F\}$.
\end{defn}

\begin{proposition}\label{lowerL}
For any graph\/ $F$, if a\/ $q$-graph\/ $H$ contains a\/ $(q+1)$-copy of\/ $F$, then\/ $S_{H^L}$ must contain a copy of\/ $F_f$ for some 1-selection\/ $f$ of\/ $F$.
\end{proposition}

\begin{proof}
If $H$ contains a $(q+1)$-copy $H'$ of $F$ with support $V$, then for any $v\in V$ there can be at most one edge $e\in E(F)$ with $v\in e$ such that the $q$-edge $\x$ playing the role of $e$ in $H'$ has $x_v<\frac{q+1}{2}$. (Otherwise the two $q$-edges would not intersect at $v$.) So we can define the 1-selection $f$ of $F$ by letting $f(v)$ this edge. Removing the corresponding $q$-edges leaves a $(q+1)$-copy of $F_f$ in $H^L$.  
\end{proof}

\begin{lemma}\label{chi1chi}
For any graph\/ $F$ containing a component\/ $C$ with\/ $|V(C)|<|E(C)|$, we have\/ $$\ex(n,F,q)\ge \left[q^2-(q-\lceil \frac{q+1}{2}\rceil +1)^2\frac{1}{\chi_1(F)-1}+o(1)\right]\binom{n}{2}.$$ 
In particular, if\/ $\chi(F)=\chi_1(F)$, then\/ $$\ex(n,F,q)=  \left[q^2-(q-\lceil \frac{q+1}{2}\rceil +1)^2\frac{1}{\chi(F)-1}+o(1)\right]\binom{n}{2}$$ holds.
\end{lemma}

\begin{proof}
To see the lower bound consider the $q$-graph $H$ that we obtain by setting $S_{H^L}=F$, where $G$ is $\cF_F$-free of size $\ex(n,\cF_F)$, $H^L=\{\x:\x\in Q(n,2)^L,\ S_\x\in E(G)\}$ and $H:=(Q(n,2)\setminus Q(n,2)^L)\cup H^L$.
This $H$ cannot contain a $(q+1)$-copy of $F$ by Proposition \ref{lowerL}. The estimate on the size of $H$ follows from Theorem \ref{ESS}.

Asymptotic equality in the second part of the lemma follows from Proposition \ref{upperL} and Theorem \ref{ESS}.
\end{proof}

Lemma \ref{chi1chi} establishes the asymptotics of $\ex(n,F,q)$ for all non-unicyclic\footnote{A graph $F$ is \textit{non-unicyclic}, if it contains a component\/ $C$ with\/ $|V(C)|<|E(C)|$.} graphs $F$ with $\chi_1(F)=\chi(F)$. Unfortunately, we can only characterize the class of graphs satisfying $\chi_1(F)=\chi(F)=2$. The edge set of the image of any 1-selection forms a graph with all components $C$ having at most $V(C)$ vertices, i.e. all components are unicyclic or trees. On the other hand, for any such graph $G$, one can define the 1-selection with image exactly $E(G)$: if $x_1,x_,\dots, x_k$ is the unique cycle, then $f(x_i):=(x_ix_{i+1})$ and for vertices $y\neq x_i$, one can let $f(y):=(yz)$, where $z$ is the unique neighbor of $y$ one closer to the cycle than $y$.

\bigskip

We close this section with comparing Theorem \ref{t:tripart} to the lower bound of Lemma \ref{chi1chi} and the upper bound of Proposition \ref{upperL}
First, we determine $\chi_1(K_{r,s,t})$. A detailed introduction of and more results on the robust chromatic number $\chi_1(F)$ will be given in the forthcoming manuscript \cite{BBPTV}. Note that $\chi_1(K_{1,1,1})=\chi_1(C_3)=1$.

\begin{proposition}\label{chistatement}
For\/ $1\le r\le s\le t$ with\/ $
t\ge 2$ we have\/ $\chi_1(K_{r,s,t})=2$ if and only if\/ $r\le 2$. Otherwise\/ $\chi_1(K_{r,s,t})=\chi(K_{r,s,t})=3$.
\end{proposition}

\begin{proof}
We have $\chi_1(K_{r,s,t})\ge 2$ if and only if $t\ge 2$ as in any 1-selection there remains an edge. If $r= 2$, then let $u,v$ be the vertices in the part of size 2, and let $B$ and $C$ be the parts of size $s$ and $t$. Then for any 1-selection $f$ with $f(b)=\{u,b\}$ for all $b\in B$ and $f(c)=\{v,c\}$ for all $c\in C$, the 1-removed $(K_{r,s,t})_f$ is bipartite showing $\chi_1(K_{r,s,t})=2$ for $r\le 2$. 

For the case $r\ge 3$, it is enough to consider $G=K_{3,3,3}$. Suppose to the contrary that $\chi(G_f)\le 2$ for some 1-selection $f$. Then $G_f$ contains an independent set of size at least 5. But any set of 5 vertices induces at least 6 edges in $K_{3,3,3}$ so $f$ cannot remove all of them.
\end{proof}

Observe that for $q=2$, substituting a $\chi$- or $\chi_1$-value 2 or 3 to the lower and upper bounds of Lemma \ref{chi1chi} and Proposition \ref{upperL} would yield bounds $(3+o(1))\binom{n}{2}$ and $(\frac{7}{2}+o(1))\binom{n}{2}$. Therefore Theorem \ref{t:tripart} (2) shows that $\ex(n,K_{2,s,t},2)$  has a value midway through these bounds whenever $3\le s\le t$ holds.

\section{Proofs}

\subsection{Bipartite graphs with unicyclic components}

\begin{lemma}\label{complement}
Let\/ $G$ be a bipartite graph such that all its components are unicyclic or trees. Suppose\/ $H\subseteq Q(n,2)$ does not contain any\/ $(q+1)$-copy of graph\/ $G$. Then for any\/ $(a,b)\in [q]^2$, the graph\/ $H_{(a,b),(\overline{a},\overline{b})}$ has\/ $o(n^2)$ edges. 
\end{lemma}

\begin{proof}
The proof is based on the following simple claim.

\begin{claim}\label{compbip}
$H_{(a,b),(\overline{a},\overline{b})}$ does not contain\/ $K_{R,R}$, where\/ $R=R(|V(G)|,|V(G)|)$ is the bipartite Ramsey number, the minimum value of\/ $c$ such that in any 2-edge-coloring of\/ $K_{c,c}$ there exists a monochromatic copy of\/ $K_{|V(G)|,|V(G)|}$.  
\end{claim}

\begin{proof}[Proof of claim]
Suppose $H_{(a,b),(\overline{a},\overline{b})}$ contains a copy of $K_{R,R}$. If $a,b \ge \frac{q+1}{2}$ or $\overline{a},\overline{b} \ge \frac{q+1}{2}$, then as $G\subset K_{|V(G)|,|V(G)|}$, $H$ contains a $(q+1)$-copy of $G$. So we may and will assume that $a\le \overline{b}\le \frac{q+1}{2}\le b\le \overline{a}$ holds. 

 Color the edge $(i,j)$ of $H_{(a,b),(\overline{a},\overline{b})}$ blue if $\x,\y\in H$ with $x_i=a, \ x_j=b,\ y_i=\overline{a},\ y_j=\overline{b}$ and red otherwise, i.e.\ if $\x',\y'\in H$ with $x_j=a,\ x_i=b,\ y_j=\overline{a},\ y_i=\overline{b}$. By definition of $R$, we obtain that there exist two disjoint sets $U,V$ with $|U|=|V|=|V(G)|$ such that for any $u\in U$ and $v\in V$ the $q$-edges $\x,\y$ with $x_u=a,\ y_u=\overline{a},\ x_v=b,\ y_v=\overline{b}$
are present in $H$.
 We claim that a $(q+1)$-copy of $G$ with support in $U\cup V$ can be defined. By induction on the number of vertices in $G$, it is enough to embed one component. Let $2\ell$ be the length of the unique cycle in the component (if there exists one).
 We fix $u_1,u_2,\dots,u_\ell\in U$ and $v_1,v_2,\dots,v_\ell\in V$.
The $q$-edges $\x^j$ with $x^j_{u_j}=a$, $x^j_{v_j}=b$ and $\y^j$ with $y^j_{v_j}=\overline{b}$ and $y^j_{u_{j+1}}=\overline{a}$ form a $(q+1)$-copy of $C_{2\ell}$, where $u_{\ell+1}=u_1$. To embed the trees hanging from the cycle (or the one tree if there is no cycle at all), one can proceed greedily according to the distance to the cycle by using $q$-edges from $\overrightarrow{H}_{a,b}$ if the already embedded vertex is in $V$ while using $q$-edges from $\overrightarrow{H}_{\overline{a},\overline{b}}$ if the already embedded vertex is in $U$. As $H$ should not contain a $(q+1)$-copy of $G$, this contradiction proves the claim.
\end{proof}

The statement of the lemma now follows from Claim \ref{compbip} and Theorem \ref{kst}.
\end{proof}

\begin{proof}[Proof of Theorem \ref{bipuni}]
The lower bound is given by the construction $U_{q,n}$.

For the upper bound, let $G$ be a bipartite graph with all its component containing at most one cycle and suppose $H\subseteq Q(n,2)$ does not contain any $(q+1)$-copies of $G$. Then we partition the edges of $H$ according to their entries. If $q+1$ is even, then Lemma \ref{complement} implies that $H_{\frac{q+1}{2},\frac{q+1}{2}}$ has $o(n^2)$ edges. By Lemma \ref{complement}, we have $|H_{(a,b),(\overline{a},\overline{b})}|=o(n^2)$, and so $|\{(i,j): 1\le i<j\le n: \x \in H ~\text{with}\ x_i=a,\ x_j=b\}|+|\{(i,j): 1\le i<j\le n: \y \in H ~\text{with}\ y_i=\overline{a},\ y_j=\overline{b}\}|\le (1+o(1))\binom{n}{2}$. Summing it up over all $(a,b),(\overline{a},\overline{b})$ we obtain $|H|\le (\lfloor \frac{q^2}{2}\rfloor+o(1))\binom{n}{2}$.
\end{proof}

\subsection{Forests}

\begin{lemma}\label{treelemma}
Let\/ $T$ be a tree of radius\/ $r$. If\/ $H \subseteq Q(n,2)$ does not contain any\/ $(q+1)$-copy of\/ $T$, then for any\/ $1\le a<\frac{q+1}{2}\le b\le q$ we have\/ $|\overrightarrow{H}_{a,b}|+|\overrightarrow{H}_{\overline{a},\overline{b}}|\le (2+o(1))\binom{n}{2}$. Furthermore, if\/ $a+b=q+1$ and thus\/ $a=\overline{b}$, then\/ $|\overrightarrow{H}_{a,b}|\le t_{n,r}+o(n^2)$.
\end{lemma}

\begin{proof}

By Lemma \ref{complement}, we have $|H_{(a,b),(\overline{a},\overline{b})}|=o(n^2)$. This implies $|\overrightarrow{H}_{a,b}|+|\overrightarrow{H}_{\overline{a},\overline{b}}|\le (2+o(1))\binom{n}{2}$ if $a+b\neq q+1$.

It remains to show $|\overrightarrow{H}_{a,b}|\le t_{n,r}+o(n^2)$, if $a+b=q+1$. Let us consider $\overrightarrow{H}'_{a,b}$ that we obtain from $\overrightarrow{H}_{a,b}$ by removing directed cycles of length 2, i.e.\ removing both orientations of the edges of $H_{(a,b),(b,a)}$. As $b=\overline{a}$, $a=\overline{b}$, we have $H_{(a,b),(b,a)}=H_{(a,b),(\overline{a},\overline{b}}$, and thus the number of removed arcs is $o(n^2)$. Construct the vertex sets $V_1,V_2,\dots,V_t$ as follows: $V_1$ is the set of vertices that have in-degree less than $|T|$ in $\overrightarrow{H}'_{a,b}$, and inductively $V_{i+1}$ is the set of vertices that have in-degree less than $|T|$ in $\overrightarrow{H}'_{a,b}[V\setminus \cup_{j=1}^iV_i]$. Observe that as long as $V\setminus \cup_{j=1}^iV_j\neq \emptyset$, then $V_{i+1}$ cannot be empty. Indeed, otherwise $\overrightarrow{T}$ (a directed version of $T$ with all edges oriented towards the center) can be embedded greedily to $\overrightarrow{H}'_{a,b}[V\setminus \cup_{j=1}^iV_j]$ and thus $H$ would contain a $(q+1)$-copy of $T$. Also, if $t$ is the first index for which $V= \cup_{j=1}^tV_j$, then $t\le r$. Indeed, if $v\in V_{r+1}$, then $\overrightarrow{T}$ can be embedded greedily into $\{v\}\cup \bigl( \cup_{j=1}^rV_j \bigr) $ with $v$ playing the role of the center.

We can now bound $|\overrightarrow{H}_{(a,b)}|$. By definition, there are at most $(|T|-1)n$ arcs in $\cup_{j=1}^t\overrightarrow{H}'_{a,b}[V_j]$. Let $\overrightarrow{H}''_{a,b}$ be the graph obtained by removing these arcs. As shown above, $\overrightarrow{H}''_{a,b}$ is $r$-partite, so $|\overrightarrow{H}''_{a,b}|\le t_{n,r}$. By our bounds on the number of removed arcs, we have $|\overrightarrow{H}_{a,b}|\le |\overrightarrow{H}''_{a,b}|+o(n^2)\le t_{n,r}+o(n^2)$ as claimed.
\end{proof}

\begin{proof}[Proof of Theorem \ref{gentree}] Fix $T$ and $q$, and let $H\subseteq Q(n,2)$ be a $q$-graph not containing any $(q+1)$-copies of $T$. As $S_{H^L}$ cannot contain $T$, we must have $|H^L|=O_{q,T}(n)$. We do not bound $H^S$ and assume $Q(n,2)^S\subseteq H$. In order to prove the theorem, we group the possible weights: $(a,b)$ with $(b,a)$, $(\overline{a},\overline{b})$, and $(\overline{b},\overline{a})$. It is enough to show that there exists $c=c_{q,T}$ such that for any $a,b$ with $1\le a<\frac{q+1}{2}\le b\le q$ we have $|\overrightarrow{H}_{a,b}|+|\overrightarrow{H}_{\overline{a},\overline{b}}|\le (2+o(1))\binom{n}{2}$, and if $a+b=q+1$, then $|\overrightarrow{H}_{a,b}|\le (1-c)\binom{n}{2}$.
This follows by the statement of Lemma \ref{treelemma}.
\end{proof}

\begin{proof}[Proof of Theorem \ref{2tree}]
For the lower bounds, we have the following general constructions. Let $\cA=(A_1,A_2,\dots,A_r)$ be a fixed partition of $[n]$, and for any $u\in [n]$ we set $a(u)=i$ if $u\in A_i$. Then we define two families of $q$-edges $\cF_\cA,\cF'_\cA\subseteq Q(n,2)$: $$\cF_\cA:=\{\x:\sum_i x_i=2\}\cup \{\x: x_u=1,x_v=2, \ a(u)<a(v)\}$$
and
$$\cF'_\cA:=\cF_A\setminus \{\x:x_u=x_v=1, \ a(u)=a(v)=r\}.$$
This $\cF'_A$ does not contain 3-copies of any tree $T$ of radius $r$, while $\cF_A$ does not contain 3-copies of any tree $T$ of radius $r$ and diameter $2r$, as the longest path in $\cF'_\cA$ is of length $2(r-1)$, while the longest path in $\cF_\cA$ is of length $2r-1$. Indeed, if a path contains a $q$-edge with two 1-entries, then it can be only continued by paths corresponding to paths in $\overrightarrow{H}_{1,2}$, so the length is at most $2(r-1)+1$. In $\cF'_A$ there is no such $q$-edge with support completely in $A_r$ and thus the length of the longest path is indeed $2(r-1)$.

To see the upper bounds, observe first that if $H\subseteq 2(n,2)$ does not contain a $3$-copy of $T$, then the number of $q$-edges in $H$ with two 2-entries is linear as $H_{2,2}$ cannot contain $T$. Also, by Lemma \ref{complement}, $H_{(1,2),(2,1)}$ has $o(n^2)$ edges, therefore $|H|\le |H_{1,1}|+|\overrightarrow{H}_{1,2}|+o(n^2)$. Lemma \ref{treelemma} yields the upper bound of (1). To obtain the upper bound of (2), assume that the diameter of $T$ is $2r-1$. We use the partition of $\overrightarrow{H}_{1,2}$ as in Lemma \ref{treelemma}: $V_1$ is the set of those vertices that have in-degree at most $|T|$ in $\overrightarrow{H}_{1,2}$ and $V_{i+1}$ is the set of vertices in $[n]\setminus \cup_{j=1}^iV_j$ that have in-degree at most $|T|$ in $\overrightarrow{H}_{1,2}[[n]\setminus \cup_{j=1}^iV_j]$. Again, as long as $[n]\setminus \cup_{j=1}^iV_j\neq \emptyset$, the set $V_{i+1}$ is non-empty, otherwise we would find a 3-copy of $T$. Also, $[n]=\cup_{i=1}^rV_r$ as otherwise we would find a 3-copy of $T$ in $H$. Finally, $H_{1,1}[V_r]$ has maximum degree less than $|T|$ as otherwise again one would be able to greedily embed $T$ in $H$ using the vertex $v\in V_r$ with degree at least $|T|$ and its neighbors in $H_{1,1}[V_r]$ and directed trees in $\overrightarrow{H}_{1,2}$ rooted at the neighbors.

We obtained that apart from a $o(n^2)$ error term $|\overrightarrow{H}_{1,2}|$ is bounded by the number of edges in a complete $r$-partite graph $G$, and $|H_{1,1}|$ is bounded by the number of pairs not completely inside the last part $V_r$ of $G$. If one fixes the size of $V_r$, then standard symmetrization shows that maximum is obtained when all other parts have (almost) equal size. So we have $r-1$ parts of size $\alpha n$ and one part of size $(1-(r-1)\alpha)n$. Then apart from a sub-quadratic error term, we have
$$|H_{1,1}|+|\overrightarrow{H}_{1,2}|\le \binom{n}{2}-\binom{(1-(r-1)\alpha)n}{2}+\binom{n}{2}-(r-1)\binom{\alpha n}{2}-\binom{(1-(r-1)\alpha)n}{2}.$$
It is easy to see that the above expression takes its maximum at $\alpha=\frac{2n}{2r-1}$ proving (2).
\end{proof}

\subsection{Odd cycles}

In this subsection, we prove Theorem \ref{triangle}.  Let us remind the reader that all lower bounds in this subsection follow from $\ex(n,F,q)\ge q^2\cdot \ex(n,F)$. By Mantel's theorem, we have $\ex(n,C_3)=\lfloor \frac{n^2}{4}\rfloor$. To obtain the upper bound, we start with some general lemmas.

\begin{lemma}\label{atlon}
Suppose a\/ $q$-graph\/ $H \subset Q(n,2)$  does not contain a\/ $(q+1)$-copy of\/ $C_3$. Let\/ $1\le a\le b\le q$ be integers with\/ $a+b=q+1$.

\medskip 

\begin{enumerate}
    \item 
    If\/ $a=b=\frac{q+1}{2}$, then\/ $$|H_{\frac{q+1}{2},\frac{q+1}{2}}|\le \lfloor \frac{n^2}{4}\rfloor.$$
    \item
    If\/ $a<b$, then\/ $$|\overrightarrow{H}_{a,b}|\le 2\cdot\lfloor \frac{n^2}{4}\rfloor.$$
\end{enumerate}
\end{lemma}

\begin{proof}
Mantel's theorem and the fact that $H_{\frac{q+1}{2},\frac{q+1}{2}}$ is triangle-free imply (1).

To see (2), we proceed by induction on $n$ with the base cases $n=2,3$ being left to the reader. If $H_{(a,b),(b,a)}$ is empty, then $|\overrightarrow{H}_{a,b}|\le \binom{n}{2}\le 2\cdot \lfloor \frac{n^2}{4}\rfloor$. If $(ij)\in H_{(a,b),(b,a)}$, then for any $u\neq i,j$ at most two of $(iu),(ui),(ju),(uj)$ belong to $\overrightarrow{H}_{a,b}$. Indeed, $(iu)\in \overrightarrow{H}_{a,b}$ implies $(uj)\notin \overrightarrow{H}_{a,b}$ as then the oriented 3-cycle $(iu),(uj),(ji)$ in $\overrightarrow{H}_{a,b}$ would correspond to a $(q+1)$-copy of a triangle in $H$. An analogous statement shows that $(ui)\in \overrightarrow{H}_{a,b}$ implies $(ju)\notin \overrightarrow{H}_{a,b}$. Applying induction to $\{\x\in H: i,j\notin S_\x\}$, we obtain 
\[
|H|\le 2+2(n-2)+2\cdot \lfloor \frac{(n-2)^2}{4}\rfloor=2\cdot\lfloor \frac{n^2}{4}\rfloor.
\]
\end{proof}

\begin{lemma}\label{atlonkivul}
Suppose a\/ $q$-graph\/ $H \subset Q(n,2)$  does not contain a\/ $(q+1)$-copy of a\/ $C_3$ and\/ $\frac{q+1}{2}< a\le q$.
Then\/ $$|H_{a,a}|+|H_{\overline{a},\overline{a}}|\le 2\cdot\lfloor \frac{n^2}{4}\rfloor.$$
\end{lemma}

\begin{proof}
The argument is very similar to the proof of Lemma \ref{atlon} (2). We proceed by induction on $n$ and we leave the base cases $n=2,3$ to the reader. Observe that if $H_{a,a}$ is empty, then $|H_{a,a}|+|H_{\overline{a},\overline{a}}|\le \binom{n}{2}\le 2\cdot \lfloor \frac{n^2}{4}\rfloor$. If $\{i,j\}\in H_{a,a}$ (and so for $\x\in Q(n,2)$ with $x_i=x_j=a$ we have $\x\in H$), then for any $u\neq i,j$ $\{u,i\}\in H_{\overline{a},\overline{a}}$ ($\{u,j\}\in H_{\overline{a},\overline{a}}$) implies $\{u,j\}\notin H_{a,a}$ ($\{u,i\}\notin H_{a,a}$) as the $q$-edges corresponding to these containments together with $\x$ would form a $(q+1)$-copy of a triangle. Thus, by induction, we obtain 
\[
|H_{a,a}|+|H_{\overline{a},\overline{a}}|\le 2+2(n-2)+2\cdot \lfloor \frac{(n-2)^2}{4}\rfloor=2\cdot \lfloor \frac{n^2}{4}\rfloor.
\]
\end{proof}

Now we are ready to prove Theorem \ref{triangle} on $\ex(n,C_3,2)$.

\begin{proof}[Proof of Theorem \ref{triangle}]
Suppose $H\subseteq 2(n,2)$ does not contain a 3-copy of a triangle. By applying Lemma \ref{atlon} (2) to bound $|\overrightarrow{H}_{1,2}|$ and Lemma \ref{atlonkivul} to bound $|H_{1,1}|+|H_{2,2}|$, we obtain $|H|=|\overrightarrow{H}_{1,2}|+|H_{1,1}|+|H_{2,2}|\le 4\cdot \lfloor \frac{n^2}{4}\rfloor$.
\end{proof}

Let us make some comments on what makes the problem of determining $\ex(n,C_3,q)$ more difficult for larger values of $q$. If $q=3$, then Lemma \ref{atlon} (1) and (2) yield $|H_{2,2}|\le \lfloor \frac{n^2}{4}\rfloor$ and $|\overrightarrow{H}_{1,3}|\le 2\cdot \lfloor \frac{n^2}{4}\rfloor$. By applying Lemma \ref{atlonkivul}, one obtains $|H_{1,1}|+|H_{3,3}|\le 2\cdot \lfloor \frac{n^2}{4}\rfloor$. It remains to bound $|\overrightarrow{H}_{1,2}|$ and $|\overrightarrow{H}_{2,3}|$.
As $H_{a,b}$ must be $C_3$-free for any $\frac{q+1}{2}\le a\le b\le q$, the size of $\overrightarrow{H}_{2,3}$ is clearly at most $2\cdot\lfloor \frac{n^2}{4}\rfloor$, but one cannot derive any non-trivial bound on $|\overrightarrow{H}_{1,2}|$ without further conditions. Thus, one might try to bound $|\overrightarrow{H}_{1,2}|+|\overrightarrow{H}_{2,3}|$ and in general, $|\overrightarrow{H}_{a,b}|+|\overrightarrow{H}_{\overline{a},\overline{b}}|$ and show that this sum is at most $4\cdot \lfloor \frac{n^2}{4}\rfloor$. Unfortunately, this is not necessarily true even if $\frac{q+1}{2}\le a<b$. Namely, if $A\subset [n]$ with $|A|=3n/4$ and $B=[n]\setminus A$, then one can set $$\overrightarrow{H}_{a,b}=\{(ij):i\in A, j\in B\}, \  \overrightarrow{H}_{\overline{a},\overline{b}}=\{(ij):|\{i,j\}\cap B|\le 1\}.$$ 
It is easy to see that the corresponding $q$-graph does not contain a $(q+1)$-copy of $C_3$ and $|\overrightarrow{H}_{a,b}|=\frac{3n^2}{4}$, $|\overrightarrow{H}_{\overline{a},\overline{b}}|=2(\binom{n}{2}-\binom{n/4}{2})$, so $|\overrightarrow{H}_{a,b}|+|\overrightarrow{H}_{\overline{a},\overline{b}}|= (\frac{9}{8}+o(1))n^2>4\lfloor \frac{n^2}{4}\rfloor$. We were able to prove that $(\frac{9}{8}+o(1))n^2$ is an upper bound, and with this we have $9\lfloor \frac{n^2}{4}\rfloor\le \ex(n,C_3,3)\le (9.5+o(1))\lfloor \frac{n^2}{4}\rfloor$.

Things get even worse for $q\ge 4$ as then pairs $a,b$ appear with $a+b>q+1$ and $a<\frac{q+1}{2}$ and for these values of $a$ and $b$, we do not have even the $(\frac{9}{8}+o(1))n^2$ upper bound on $|\overrightarrow{H}_{a,b}|+|\overrightarrow{H}_{\overline{a},\overline{b}}|$.

\subsection{Graphs having a component with more edges than vertices}

We start with the proof of Theorem \ref{random}. We will need a result from the book \cite{JLR} of Janson, \L uczak, and Ruci\'nski. The setting is as follows: let $\Gamma$ be a finite set (in our case the edge set of $K_{m,m,\dots,m}$), and we have a set of independent Bernoulli random variables corresponding to the elements $\gamma\in \Gamma$. Then for any subset $A\subset \Gamma$, we write $I_A$ to denote the indicator variable of the event that for every $\gamma \in A$ the corresponding Bernoulli random variable takes value 1. (For us $A$ will be the edge set of an $r$-clique in $K_{m,m,\dots,m}$.) For some $\cS\subseteq 2^{\Gamma}$ we define the variable $X=\sum_{A\in \cS}I_A$. (We will use $\cS$ as the family of edge sets of all $r$-cliques in $K_{m,m,\dots,m}$.)

\begin{theorem}[Theorem 2.14 in \cite{JLR}]\label{jlr}
Let\/ $X=\sum_AI_A$ and\/ $\overline{\Delta}=\sum\sum_{A\cap B\neq \emptyset}\mathbb{E}(I_AI_B)$. Then for\/ $0\le t \le \mathbb{E}X$
\[
\mathbb{P}(X\le \mathbb{E}X-t)\le \exp \left(-\frac{t^2}{2\overline{\Delta}}\right).
\]
\end{theorem}

\begin{proof}[Proof of Theorem \ref{random}]
An edge can be contained in at most $m^{r-2}$ $r$-cliques in $K(m,r,p)$. Thus the number of $r$-cliques that a 1-selection can remove from $K(m,r,p)$ is at most $rm^{r-1}$. We apply Theorem \ref{jlr}, with $X$ being the number of $r$-cliques in $K(m,r,p)$. Clearly, we have  $\mathbb{E}X=m^rp^{\binom{r}{2}}=\omega(m^{r-1})$ by the assumption $p=\omega(m^{-1/\binom{r}{2}})$. We set $t=\frac{\mathbb{E}X}{2}$.
If the edge set of two cliques intersect, then they share at least two vertices. Therefore, we have
\[
\overline{\Delta}=\sum_{i=2}^r\binom{r}{i}m^i[m(m-1)]^{r-i}p^{2\binom{r}{2}-\binom{i}{2}}.
\]
Writing $\beta_i=\binom{r}{i}m^i[m(m-1)]^{r-i}p^{2\binom{r}{2}-\binom{i}{2}}$, we have $\frac{\beta_{i+1}}{\beta_i}=\frac{r-i}{i+1}\cdot\frac{m}{m(m-1)}\cdot p^{-2i}$, which is monotone increasing in $i$. Therefore $\beta_i$ is maximized either at $i=2$ or at $i=r$. Clearly, $\beta_r=2\mathbb{E}X$ and $\beta_2=(1+o(1))\binom{r}{2}\cdot \frac{(\mathbb{E}X)^2}{m^2p}$. So as long as $\mathbb{E}X$ and $m^2p$ tend to infinity, so does $\min\{\frac{t^2}{r\beta_2},\frac{t^2}{r\beta_r}\}$ which bounds $\frac{t^2}{2\overline{\Delta}}$ from below. As the assumption $p=\omega(m^{-1/\binom{r}{2}})$ implies both, applying Theorem \ref{jlr} we obtain that with high probability the number of $r$-cliques in $K(m,r,p)$ is greater than $rm^{r-1}$ and thus $\chi(K(m,r,p)_f)\ge r$ for any 1-selection $f$.
\end{proof}

\bigskip

Now we turn to Theorem \ref{t:tripart}, the three parts of which will be proved separately. Proposition \ref{chistatement} states that if $3\le r\le s\le t$, then $\chi_1(K_{r,s,t})=\chi(K_{r,s,t})=3$ holds. In this case, Lemma \ref{chi1chi} implies the following result which is more general than part (3) of Theorem \ref{t:tripart}.

\begin{corollary}
For any\/ $3\le r\le s\le t$, we have\/ $$\ex(n,K_{r,s,t},q)= \left[q^2-(q-\lceil \frac{q+1}{2}\rceil +1)^2\frac{1}{2}+o(1)\right]\binom{n}{2}.$$ 
\end{corollary}

The next lemma establishes the first half of part (1) of Theorem \ref{t:tripart}

\begin{lemma}
For any\/ $1\le s$ and\/ $2\le t$ we have\/ $\ex(n,K_{1,s,t},2)=(3+o(1))\binom{n}{2}$.
\end{lemma}

\begin{proof}
We fix $s,t$, we write $f(n)=\ex(n,K_{1,s,t},2)$ and for a fixed $\varepsilon>0$ we want to show $f(n)\le (3+\varepsilon)\binom{n}{2}$. 
We fix $T$ such that $\frac{t}{T}\le \varepsilon/2$. Suppose $H\subseteq 2(n,2)$ does not contain any 3-copies of $K_{1,s,t}$ and $|H|=f(n)$. If $e(H_{2,2})\le \varepsilon \binom{n}{2}$, then clearly $|H|\le (3+\varepsilon)\binom{n}{2}$. Otherwise, by Theorem \ref{kst}, $H_{2,2}$ contains a copy $K$ of $K_{T,T}$ with vertex set $U=U_1\cup U_2$ with the $U_i$s forming the two parts of $K$. Observe that for any $v\notin U$, we have either $|N^-_{\overrightarrow{H}_{1,2}}(v)\cap U_1|<t$ or $|N^-_{\overrightarrow{H}_{1,2}}(v)\cap U_2|<t$, as otherwise $v$ and the incident $(1,2)$-edges together with $K$ would form a 3-copy of $K_{1,s,t}$. Therefore, for any $v\notin U$, the number of edges $\x$ in $H$ with $x_v\neq 0$ and $S_\x\cap U\neq \emptyset$ is at most $4(T+t-1)+2(T-t+1)$. Considering edges of $H$ within $U$, between $U$ and $V(H)\setminus U$, and outside $U$, yields
\[
f(n)\le 4\binom{2T}{2}+[3\cdot (2T)+2(t-1)](n-2T)+f(n-2T)\le 4\binom{2T}{2}+(3+\varepsilon/2)2T(n-2T)+f(n-2T).
\]
Applying this repeatedly, we obtain
\[
f(n)\le 4\binom{2T}{2}\frac{n}{2T}+(3+\varepsilon/2)2T\frac{n(\frac{n}{2T}-1)}{2}\le (3+\varepsilon)\binom{n}{2},
\]
if $n$ is large enough.
\end{proof}

The next result is the second part of part (1) of Theorem \ref{t:tripart}.

\begin{theorem}\label{22t}
For any\/ $1\le r\le s\le 2\le t$ we have\/ $\ex(n,K_{2,s,t},2)=(3+o(1))\binom{n}{2}$.
\end{theorem}

\begin{proof}
It suffices to consider $r=s=2$ and unrestricted large $t$.
In the argument $\varepsilon_1,\varepsilon_2,\dots$ will denote positive constants, where $\varepsilon_1$ is arbitrarily small and the others are chosen appropriately.

We apply induction on $n$.
As an anchor, up to any fixed threshold $n_0$ it can be assumed that $\ex(n,K_{2,s,t},2)\le 3\binom{n}{2} + C$ for a suitably chosen constant $C$, for all $n\le n_0$.

Let now $n>n_0$ and let $H\subseteq 2(n,2)$ be any 2-graph with $e(H)\ge (3+ \varepsilon_1) \binom{n}{2}$.
Then clearly $e(H_{2,2})\ge \varepsilon_1 \binom{n}{2}$ holds, and as above, $H_{2,2}$ contains a copy of $K_{T,T}$ with partite classes $U_1$ and $U_2$ of size $T$, where $T$ can be chosen huge, e.g.\ exponential in both $t$ and $1/\varepsilon_1$.

If $H_{2,2}$ has at most $\varepsilon_2 n$ edges from $U:=U_1\cup U_2$ to its complement, then induction proves that $H$ contains a 3-copy of $K_{2,2,t}$.
Otherwise consider the set $Y$ of vertices $y\notin U$ that have at least $\varepsilon_3 |U|$ neighbors in $U$ within $H_{2,2}$.
We may assume that at least half of $Y$ have at least $\varepsilon_3 |U_1|$ $(2,2)$-neighbors in $U_1$.
Then there exist a $U'\subset U$ of size $t$ and a $Y'\subset Y$ of size $\varepsilon_4 n$ completely joined in $H_{2,2}$.
If $\varepsilon_5 n$ of those vertices of $Y'$ have more than one $(2,2)$-neighbor in $U_2$, then some vertex pair of $U_2$ has many common neighbors in $Y'$, hence $H$ contains a 3-copy of $K_{2,2,t}$. (In fact, even a $K_{2,t,t}$ occurs.)
Else there is a subset $Y''\subset Y'$ of $\varepsilon_5 n$ vertices having at most one $(2,2)$-neighbor in $U_2$.

Consider now the bipartite graph $B=H_{(1,2),(2,1)}[Y'',U_2]$. 
By the choice of $Y$, each $y\in Y''$ has degree at least $2\varepsilon_3 |U_2|$ in $B$, because there cannot be more than $4|U_1|$ edges from $y$ to $U_1$ in $H$.
Thus $B$ contains a $C_4=y'u'y''u''$, and consequently a 3-copy of $K_{2,2,t}$ with vertex classes $\{y',y''\}$, $\{u',u''\}$ and $U'$ occurs in $H$.
\end{proof}

Let us state the second part of Theorem \ref{t:tripart} separately.

\begin{theorem}\label{2st}
For any\/ $3\le s\le t$ we have\/ $\ex(n,K_{2,s,t},2)=(\frac{13}{4}+o(1))\binom{n}{2}$.
\end{theorem}

Before the proof of Theorem \ref{2st}, we need some preparations.

\begin{defn}
We say that\/ $e$ is an\/ $i$-edge of\/ $W_k:E(K_k)\rightarrow \{0,2,3\}$ if\/ $W_k(e)=i$. A\/ $W_k$ satisfies \textit{Condition\/ $(\star)$} if\/
\begin{itemize}
\item[(i)] the graph of 3-edges is triangle-free and
\item[(ii)] if\/ $uv,vx,xy,yu$ 3-edges form a\/ $C_4$, then at least one of\/ $ux$ and $vy$ is a 0-edge.
\end{itemize}
The total weight $w(W_k)$ of $W_k$ is $\sum_{e\in E(K_k)}W_k(e)$. The degree of $v$ in $W$ is $d(v)=d_W(v)=\sum_{x\neq v}W_k(vx)$.
\end{defn}

\begin{lemma}\label{023}
If\/ $W_k:E(K_k)\rightarrow \{0,2,3\}$ satisfies\/ $(\star)$, then\/ $w(W_k)\le (\frac{9}{4}+o(1))\binom{k}{2}$, and this bound is asymptotically sharp.
\end{lemma}

\begin{proof}
We start with an observation that applies Zykov's symmetrization in this setting.

\begin{claim}\label{zykov}
If\/ $W$ satisfies\/ $(\star)$, then, for any pair\/ $u,v$ of vertices, so does\/ $W^{u\rightarrow v}$ defined as\/ $W^{u\rightarrow v}(uv)=0$, $W^{u\rightarrow v}(ux)=W(vx)$ for all\/ $x\neq u,v$ and\/ $W^{u\rightarrow v}(e)=W(e)$ for all other edges\/ $e\in E(K_k)$.

Furthermore, if\/ $uv$ is an\/ $i$-edge and\/ $d_W(v)-d_W(u)\ge i$, then\/ $w(W)\le w(W^{u\rightarrow v})$ and if\/ $\/d_W(v)-d_W(u)> i$, then\/ $w(W)> w(W^{u\rightarrow v})$.
\end{claim}

\begin{proof}[Proof of Claim]
Clearly, $W^{u\rightarrow v}$ satisfies Condition $(\star)$ on $V(K_k)\setminus \{u\}$ as does $W$. For $x,y\neq v$ we cannot have $W^{u\rightarrow v}(ux)=W^{u\rightarrow v}(uy)=W^{u\rightarrow v}(xy)=3$ as then we would have $W(vx)=W(vy)=W(xy)=3$ contradicting the triangle-free property of 3-edges of $W$. A similar argument works for $C_4$s in 3-edges in $W^{u\rightarrow v}$ containing $u$ but not $v$. Finally, $W^{u\rightarrow v}(uv)=0$ settles the cases of possible triangles and $C_4$ containing both $u$ and $v$.

The second part of the claim follows as $w(W^{u \rightarrow v})-w(W)=d_W(v)-d_W(u)-W(uv)$.
\end{proof}

Suppose $w=W_k$ satisfies $(\star)$. Let $v$ be a vertex with maximum degree. By repeatedly applying Claim \ref{zykov} with $W^{x\rightarrow v}$ for all $x$ with $W(vx)=0$, we obtain $W^1$ and $U_1\ni v$ such that $w(W_1)\ge w(W)$ and for any $u,u'\in U_1$ and $x\notin U_1$ we have $W^1(ux)=W^1(u'x)\ge 2$. Then we repeat this process for a vertex $v_2$ of maximum degree $d_{W^1}(v_2)$ among vertices not in $U_1$. After a finite number of such rounds we obtain $W^\ell$ and $U_1,U_2,\dots,U_\ell$ such that $w(W^\ell)\ge w(W)$ and for any $x,y,\in U_i$, $u,v\in U_j$ ($i=1,2,\dots,\ell$, $j=1,2,\dots,\ell$, $i\neq j$) we have $W^\ell(xy)=0$  and $W^\ell(xu)=W^\ell(yv)\ge 2$. In particular, $d_{W^\ell}(x)=d_{W^\ell}(y)$.

Finally, as long as there exists an $i$-edge $uv$ for $i=2,3$ with $d(v)-d(u)>i$, we again apply Claim \ref{zykov}. Every such symmetrization moves vertex $u \in U_j$ to the class $U_{j'}$ of $v$ and strictly increases the weight. Since the weight cannot be larger than $3\binom{n}{2}$, we obtain $W^*$ and partition $U^*_1,U^*_2,\dots,U^*_{\ell^*}$ of $V(K_n)$ such that all above properties hold, and also $|d_{W^*}(u)-d_{W^*}(v)|\le 3$ for any pair $u,v$ of vertices. This implies that either $w(W)\le w(W^*)\le (\frac{9}{4}+o(1))\binom{k}{2}$ or $d_{W^*}(v)>\frac{9}{4}k$ holds for all $v$.

Without loss of generality, we may assume that $U^*_{\ell^*}$ is of largest size, which we denote by $x$, among all classes, and that $U^*_1,U^*_2,\dots, U^*_m$ are those other classes that are connected by 3-edges to $U^*_{\ell^*}$ and other classes are connected by 2-edges.  Let $y_1,y_2,\dots,y_m$ denote the sizes of $U^*_1,U^*_2,\dots, U^*_m$. Observe the following.

\medskip 

\begin{itemize}
    \item 
    As we only have 0-edges within $U^*_{\ell^*}$ and degrees are all at least $\frac{9}{4}k$, we must have $\sum_{i=1}^my_i\ge \frac{k}{4}+2x$.
    \item
    The triangle-free part of $(\star)$ implies that all edges between $U_i$ and $U_j$ ($1\le i<j\le m$) are 2-edges.
    \item
    The $C_4$ part of $(\star)$ implies that if $z$ is connected by 3-edges to both $u_i\in U^*_i$ and $u_j\in U^*_j$ ($1\le i<j\le m$), then $z\in U^*_{\ell^*}$.
    \item
    Therefore, if $V_i$ denotes $\{z\notin U^*_{\ell^*}:W^*(zu_i)=3 ~\forall u_i\in U^*_i\}$ for $i=1,2,\dots,m$, then $U^*_{\ell^*},U^*_1,U^*_2,\dots,U^*_m,V_1,V_2,\dots,V_m$ are pairwise disjoint.
    \item
    As all degrees are at least $\frac{9}{4}k$, we must have $|V_i|>\frac{k}{4}+2y_i-x$.
\end{itemize}

\medskip 

Because of the pairwise disjoint property, we can add up all lower bounds on sizes and obtain
\[
k\ge x+\sum_{i=1}^my_i +\sum_{i=1}^m(\frac{k}{4}+2y_i-x)=\frac{mk}{4}+3\sum_{i=1}^my_i-(m-1)x>\frac{(m+3)k}{4}-(m-7)x.
\]
The right-hand side is at least $k$ (which is a contradiction) unless $m=0$ or $x>\frac{k}{4}$. But if $m=0$, then the degrees in $U^*_{\ell^*}$ are at most $2\binom{k}{2}$, while if $x>\frac{k}{4}$, then $|\cup_{i=1}^mU^*_i|> \frac{k}{4}+2x>\frac{3k}{4}$. This is impossible as $U^*_{\ell^*}$ and $\cup_{i=1}^mU^*_i$ are disjoint. These contradictions finish the proof of the lemma.
\end{proof}

\begin{defn}
In a graph $G$, for disjoint subsets $X,Y$ of $V(G)$ we write $d(X,Y)$ to denote the number of edges $e$ in $E(G)$ with $|e\cap X|=|e\cap Y)|=1$. For\/ $\varepsilon > 0$ a pair\/ $X,Y$ of disjoint subsets of vertices in\/ $G$ is\/ \textit{$\varepsilon$-regular} if for\/ $A\subseteq X$, $B\subseteq Y$ with\/ $|A|\ge \varepsilon |X|$, $|B|\ge \varepsilon |Y|$ we have\/ $|d(X,Y)-d(A,B)|\le \varepsilon$.
\end{defn}

\begin{theorem}[Szemer\'edi's Regularity Lemma \cite{Sz}]\label{reg}
For any positive integer\/ $m$ and positive real\/ $\varepsilon$ there exists an integer\/ $M$ such that any graph\/ $G$ on at least\/ $M$ vertices can be partitioned into\/ $V_0,V_1,\dots,V_k$ with\/ $|V_0|\le \varepsilon n$, $|V_1|=|V_2|=\dots =|V_k|$ such that\/ $m\le k \le M$ and all but\/ $\varepsilon k^2$ pairs\/ $V_i,V_j$ are\/ $\varepsilon$-regular.  
\end{theorem}

\begin{defn}
For any graph\/ $G$ and integer\/ $t\ge 2$, we denote by\/ $G(t)$ its\/ \textit{$t$-blow-up} obtained by replacing each vertex\/ $v$ of\/ $G$ by a\/ $t$-element independent set\/ $I_v$ and vertices\/ $v' \in I_v, u'\in I_u$ are adjacent in\/ $G(t)$ if and only if\/ $u$ and\/ $v$ are adjacent in\/ $G$.
\end{defn}

\begin{lemma}[Key Lemma \cite{KS}]\label{key}
Given\/ $d>\varepsilon>0$, a graph\/ $R$, and a positive integer\/ $m$, let the graph\/ $G$ be constructed by replacing every vertex of\/ $R$ by\/ $m$ vertices, and every edge of\/ $R$ by an\/ $\varepsilon$-regular pair of density at least\/ $d$. Let\/ $H$ be a subgraph of\/ $R(t)$  with maximum degree\/ $\Delta$. Let\/ $\delta=d-\varepsilon$ and\/ $\varepsilon_0=\delta^\Delta/(2+\Delta)$. If\/ $\varepsilon\le \varepsilon_0$ and\/ $t-1\le \varepsilon m$, then\/ $H\subset G$.
\end{lemma}

\begin{proof}[Proof of Theorem \ref{2st}] 
We start with the construction $H\subset 2(n,2)$ that shows $\ex(n,K_{2,s,t},2)\ge \frac{13}{4}\binom{n}{2}$. Let $A=[\lfloor \frac{n}{4}\rfloor]$ and $B=[n]\setminus A$. For an edge $e\in \binom{A}{2}$, we set $W_H(e)=\{(1,1)\}$, for an edge $e\in \binom{B}{2}$ we set $W_H(e)=\{(1,1),(2,1),(1,2)\}$, and finally for edges $e$ with $|e\cap A|=|e\cap B|=1$ we set $W_H(e)=[2]\times [2]$. Observe that $H$ does not contain a 3-copy of $K_{2.3,3}$. Indeed, if $X$ is an 8-set of the union of the supports, then $|X\cap A|\le 3$ as $X\cap A$ can contain the support of at most one edge. Therefore $X\cap B$ must contain at least 5 vertices and thus should contain a 3-copy of $K_{2,3}$, which is impossible by Lemma \ref{nagyel}.

Suppose $H\subset 2(n,2)$ does not contain any 3-copy of $K_{2,s,s}$ with $n$ large enough.
For any positive real $d$, we want to prove $|H|\le (\frac{13}{4}+d)\binom{n}{2}$. Trivially, $|H_{1,1}|\le \binom{n}{2}$, so it is enough to show that $2|H_{2,1}|+|H_{2,2}|\le (\frac{9}{4}+d)\binom{n}{2}$ holds.

Fix $\varepsilon_1<d/20$ so small that with $\delta=d/20-\varepsilon_1$ and $\Delta=3$ the inequality $\varepsilon_1\le \delta^\Delta/(2+\Delta)$ holds.  We let $m\ge\lceil \frac{100}{d}\rceil$ and so large that Lemma \ref{023} implies $w(W_m)\le (\frac{9}{4}+d/2)\binom{m}{2}$ in its statement. We denote by $M_1$ the value obtained from Theorem \ref{reg} with $m$ and $\varepsilon_1$. We define $\varepsilon_2=\varepsilon_1/2M_1$ and obtain $M_2$ by applying \ref{reg} with $\varepsilon_2$ and $m$.

Theorem \ref{reg} yields an $\varepsilon_2$-regular partition $V_0,V_1,\dots,V_k$ of $H_{2,2}$ with $m\le k\le M_2$.
\begin{itemize}
    \item 
    The number of edges incident to $V_0$ is at most $\varepsilon_2n^2$ both in $H_{2,2}$ and $H_{2,1}$.
    \item
    By choice of $m$, the total number of edges within the $V_i$s is at most $d n^2/100$ in $H_{2,1}$ and $H_{2,2}$.
    \item
    The total number of edges between non-regular pairs in $H_{2,2}$ and $H_{2,1}$ is at most $\varepsilon_2 n^2$.
    \item
    The total number of edges of $H_{2,2}$ within $\varepsilon_2$-regular pairs with density at most $d/10$ is at most $\frac{d}{10}\binom{n}{2}$.
    \item
    The same holds for edges of $H_{2,1}$ in pairs $V_i,V_j$ that are $\varepsilon_2$-regular in $H_{2,2}$ and have density at most $d/10$ in $H_{2,1}$.
\end{itemize}
The contribution to $2|H_{2,1}|+|H_{2,2}|$ counted so far is not more than $\frac{d}{2}\binom{n}{2}$. We define an auxiliary $W_k: E(K_k)\rightarrow \{0,2,3\}$ as follows: we set $W_k(ij)=3$ if the pair $V_i,V_j$ is $\varepsilon_2$-regular in $H_{2,2}$ with density at least $d/10$, we  set $W_k(ij)=2$ if the pair $V_i,V_j$ is $\varepsilon_2$-regular in $H_{2,2}$ with density smaller than $d/10$ in $H_{2,2}$ but at least $d/10$ in $H_{2,1}$, finally we set $W_k(ij)=0$ otherwise. Clearly, the so far uncounted contribution to $2|H_{2,1}|+|H_{2,2}|$ is at most $(\frac{n}{k})^2w(W_k)$. Therefore, if we can show that $W_k$ satisfies Condition $(\star)$, then, by Lemma \ref{023}, this contribution is at most $(\frac{n}{k})^2(\frac{9}{4}+d/2)\binom{k}{2}\le (\frac{9}{4}+d/2)\binom{n}{2}$ as claimed.

The 3-edges of $W_k$ form a triangle-free graph, as otherwise there would exist $V_i,V_j, V_h$ with all pairs $\varepsilon_2$-regular of density at least $d/10$ and so for large enough $n$ Lemma \ref{key} would yield a copy of $K_{s,s,s}$ in $H_{2,2}$ and the corresponding $\x$s in $H$ would form a $3$-copy of $K_{s,s,s}$ in $H$.

Suppose towards a contradiction that $V_1,V_2,V_3,V_4$ are such that the pairs $(V_1,V_2)$, $(V_2,V_3)$, $(V_3,V_4)$, and $(V_1,V_4)$ are $\varepsilon_2$-regular of density at least $d/10$ in $H_{2,2}$ (so the corresponding 3-edges of $W_k$ form a 4-cycle) and the pairs $(V_1,V_3)$ and $(V_2,V_4)$ have density at least $d/10$ in $H_{2,1}$ (so they are 2-edges of $W_k$). Then we apply Theorem \ref{reg} to the bipartite graphs $G_1:=H_{2,1}[V_1,V_3]$ and $G_2:=H_{2,1}[V_2,V_4]$ with $\varepsilon_1$ and $m$.
Denote by $U_1,U_1,\dots,U_\ell$ the obtained partition of $G_1$,
with $m\le \ell \le M_1$. As the number of edges between $\varepsilon_1$-regular pairs of density at most $d/15$ contains at most two thirds of the edges of $G_1$, while edges between non-regular pairs or within parts or adjacent to $U_0$ are even smaller in number, therefore there must exist an $\varepsilon_1$-regular pair $U_1,U_3$ with density at least $d/15$. By the definition of regular pairs, it cannot happen that $V_1$ contains more than $\varepsilon_1|U_1|$ vertices both from $U_1$ and $U_3$. Indeed, $G_1$ is bipartite, so there would be no edge between those two subparts, while the density should not differ from the density of $U_1,U_3$ by more than $\varepsilon_1$. The analogous statement holds for $V_3$. As a consequence, possibly after relabelling, we can find $U^*_1\subset U_1\cap V_1$, $U^*_3\subset U_3\cap V_3$ with $|U^*_1|=|U^*_3|=|U_1|/2\ge |V_1|/2M_1$, so the pair $U^*_1,U^*_3$ is $2\varepsilon_1$-regular in $H_{2,1}$ with density at least $d/20$. Similarly, we can find $U^*_2\subset V_2,U^*_4\subset V_4$ that form a $2\varepsilon_1$-regular pair in $H_{2,1}$ with density at least $d/20$ and $|U^*_2|=|U^*_4|=|U^*_1|$. As $|U^*_i|\ge |V_i|/2M_1$, the pairs $(U^*_i,U^*_{i+1})$ (the indices are taken modulo 4) are $2\varepsilon_2\cdot M_2=2\varepsilon_1$-regular with density at least $d/20$. 

We apply Theorem \ref{key} to the parts $U^*_1,U^*_2,U^*_3,U^*_4$ and obtain a copy $K$ of the blow-up $K_4(2s-1)$ in the graph that takes its edges from $H_{2,2}$ between $U^*_i,U^*_{i+1}$ and from $H_{2,1}$ between the other two pairs. Consider $K[U^*_1,U^*_3]=K_{2s-1,2s-1}$ and the corresponding edges of $H$. These are edges with non-zero entries 1 and 2. So there must exist a vertex $u$ that is adjacent to at least $s$ edges $\x^i$ of $H$ corresponding to edges in $K$ such that $x^i_u=2$. Without loss of generality we can assume that $u\in U_1$; then let $u^3_1,u^3_2,\dots u^3_s\in U^*_3$ denote its neighbors, i.e.\ those vertices for which $x^i_{u^3_i}=1$. Similarly, we can find $v\in U^*_2$, $u^4_1,u^4_2,\dots u^4_s\in U^*_4$ and edges $y^i$ of $H$ with $y^i_v=2$ and  $y^i_{v^4_i}=1$. But then $u$ and $v$ can form the part of size 2 of a $(q+1)$-copy of $K_{2,s,s}$ in $H$ with the $u^3_i$s and the $v^4_j$s forming the other two parts, as the remaining edges correspond to edges of $H_{2,2}$. 

We obtained a contradiction, which shows that $W_k$ must satisfy the $C_4$ part of Condition $(\star)$. This finishes the proof of the theorem. 
\end{proof}

\section{Concluding remarks and open problems}

Let us briefly summarize the problems that are left open.

\medskip

For trees $T$ (and in general for forests) we do not know the asymptotics of $\ex(n,T,q)$ if $q\ge 3$, but we suspect that the asymptotically extremal construction will be similar to what we had in the case $q=2$, i.e.\ $\overrightarrow{H}_{a,\overline{a}}$ should be a complete $r$-partite graph where $r$ is the radius of $T$.

\bigskip

The most interesting open problem would be to determine $\ex(n,C_{2k+1},q)$. We conjecture that for any $k$ and $q$ it is $q^2 \ex(n,C_{2k+1})$ at least asymptotically. Theorem \ref{triangle} verified this for $k=1$ and $q=2$. The theorem had a simple inductive proof, so let us state here that the base cases of a potential inductive proof would work.

\begin{proposition}
For\/ $n=2,3,4$ and for any\/ $q\ge 2$, we have\/ $\ex(n,K_3,q)=q^2\lfloor \frac{n^2}{4}\rfloor$.
\end{proposition}

\begin{proof}
The $n=2$ case trivially holds.

If $n=3$, the statement follows from the fact that the set of the $3 q^2$ $q$-edges in $Q(3,2)$ can be partitioned into pairwise disjoint triples, each forming a $(q+1)$-copy of a triangle, and therefore any $H\subset Q(3,2)$ without a $(q+1)$-triangle can contain at most two $q$-edges from each of those triples. If $a,b,c$ denote the three indices of the ``underlying set'', then the partition is as follows:
\begin{itemize}
    \item For any $q$ and all $j=0,1,\dots,\lfloor \frac{q-2}{3} \rfloor$ and all $i=1,2,\dots,q-1-3j$ consider
     $$ (i+2j,j+1) , \quad (q-j,i+j) , \quad (q+1-i-j,q+1-i-2j) $$
     in all the three (cyclic) positions over the triple of pairs
      $$ (a,b) , \quad (b,c) , \quad (c,a) .$$
    \item For any $q$ and all $j=0,1,\dots,\lfloor \frac{q-3}{3} \rfloor$ and all $i=1,2,\dots,q-2-3j$ consider
     $$ (i+j,i+2j+1) , \quad (q-i-2j,q-j) , \quad (j+1,q+1-i-j) $$
     in all the three (cyclic) positions over the triplet of pairs
      $$ (a,b) , \quad (b,c) , \quad (c,a) .$$
    \item If $q=3r+1$, also take
     $$ (2r+1,r+1) , \quad (2r+1,r+1) , \quad (2r+1,r+1) .$$
    \item If $q=3r+2$, also take
     $$ (r+1,2r+2) , \quad (r+1,2r+2) , \quad (r+1,2r+2) .$$
\end{itemize}

The case $n=4$ follows from the $n=3$ case by averaging. The support of every $q$-edge belongs to two triples of the indices, and there are four such triples. Therefore, we obtain $2\ex(4,K_3,q)\le 4\ex(3,K_3,q)=4\cdot 2q^2$.
\end{proof}

One potential reason for the problem not being easy is that for even values of $q$, we have found several non-isomorphic $q$-graphs, not containing $(q+1)$-copies of $K_3$, with the same asymptotic size.

\begin{itemize}
    \item
    The first construction is $U_{q,n}$, the universal tree as defined in Definition \ref{utre}.
    \item The second construction is the $K_3$-free ordinary graph with all possible $q$-edges. It means that (if we suppose that $n$ and $q$ are both even numbers) this $q$-graph has $\frac{n^2}{4}q^2$ $q$-edges.
    \item Let $n=2k$. Consider two disjoint sets $X$ and $Y$ of size $k$. In the third construction we put all $q$-edges ($a,b$) with $1 \le a,b \le \frac{q}{2}$ into $X$ and $Y$ and put all $q$-edges ($c,d$) with $1 \le c \le \frac{q}{2}$ or $1 \le d\le  \frac{q}{2}$ between $X$ and $Y$. It is easy to check that this construction also contains $\frac{n^2}{4}q^2 - O(n)$ $q$-edges.
    \item For $q=4$ we have a fourth type of construction. Let $n=3k$ and $X,Y,Z$ be three disjoint sets of size $k$. We put the following $q$-edges between the sets:
    
    \begin{itemize}
        \item inside $X$ and $Y$ we put the $q$-edges $(3,1), (2,1), (1,1), (1,2), (1,3)$;
        \item inside $Z$ we put the $q$-edges $(2,2),(2,1),(1,2),(1,1)$;
        \item between the sets $X$ and $Y$ we put all edges ($a,b$) with $1 \le a,b \le 3$;
        \item between the sets $X \cup Y$ and $Z$ we put the following
        edges: $(1,i)$ where $1 \le i \le 4$, and $(i,j)$ where $i=2,3,4$ and $j=1,2$.
    \end{itemize}
    It is easy to see that this construction also fulfils the requirements and the number of edges is $$
   9k^2 + 2\cdot 5\cdot {k \choose 2} + 10\cdot 2k^2 + 4\cdot {k \choose 2}
   =
   36k^2 - 7k
   =
   4n^2 - O(n) .
 $$
\end{itemize}

\bigskip

Lemma \ref{chi1chi} determines the order of magnitude of $\ex(n,F,q)$ whenever $\chi_1(F)=\chi(F)$, and gives lower and upper bounds if $\chi_1(F)<\chi(F)$. Theorem \ref{t:tripart} shows examples when neither of these bounds are sharp. It would be interesting to find  sufficient conditions on $F$ so that the lower or the upper bound holds tight.


\begin{thebibliography}{99}

\bibitem{BBPTV}
G. Bacs\'o, Cs. Bujt\'as, B. Patk\'os, Zs. Tuza, M. Vizer. The robust chromatic number of a graph. \textit{Manuscript}.
\bibitem{ES}
P. Erd\H os, M. Simonovits. A limit theorem in graph theory. Studia Sci. Math. Hungar. 1 (1966), 51--57.
\bibitem{ESt}
P. Erd\H os, A.H. Stone. On the structure of linear graphs. Bull. Amer. Math. Soc. 52 (1946), 1087--1091
\bibitem{FS}
Z. Füredi, M. Simonovits. The history of degenerate (bipartite) extremal graph problems. In Erdős Centennial. Springer, Berlin, Heidelberg, 169--264, 2013.
\bibitem{JLR}
S. Janson, A. Ruci\'nski, T. \L uczak.  \textit{Random graphs}. John Wiley \& Sons, 2011.
\bibitem{KS}
J. Koml\'os, M. Simonovits. Szemer\'edi's regularity lemma and its applications in graph theory. In Combinatorics, Paul Erd\H os is Eighty, Vol. 2 (Keszthely, 1993), 295--352, Bolyai Soc. Math. Stud., 2, J\'anos Bolyai Math. Soc., Budapest, 1996.
\bibitem{KST}
P. Kővári, V.T. Sós, P. Turán. On a problem of Zarankiewicz. Colloquium Mathematicum 3 (1954), 50--57.
\bibitem{PTV}
B. Patk\'os, Zs. Tuza, M. Vizer. Vector sum-intersection theorems. arxiv2305.01328 .
\bibitem{Sz}
E. Szemer\'edi. Regular partitions of graphs. In Colloques Internationaux C.N.R.S. 260, Probl\`emes Combinatoires et Th\'eorie des Graphes (Orsay, 1976), pp. 399--401.

\end{thebibliography}
\end{document}